\documentclass[a4paper,11pt]{article}
\usepackage{array}
\usepackage{theorem}
\usepackage{amsmath,amscd,amssymb}
\usepackage{latexsym}
\usepackage{epsfig}
\theorembodyfont{\sl}

\newtheorem{lemma}{Lemma}[section]

\newtheorem{proposition}[lemma]{Proposition}

\newtheorem{introthm}{Theorem}
\newtheorem{introdefinition}{Definition}
\newcommand{\DD}{\mathbb D}

\newcommand{\FF}{\mathbb F}

\newcommand{\NN}{\mathbb N}

\newcommand{\PP}{\mathbb P}
\newcommand{\QQ}{\mathbb Q}
\newcommand{\RR}{\mathbb R}

\newcommand{\ZZ}{\mathbb Z}
\newcommand{\cB}{\mathcal B}
\newcommand{\cC}{\mathcal C}
\newcommand{\cD}{\mathcal D}
\newcommand{\cE}{\mathcal E}

\newcommand{\cG}{\mathcal G}
\newcommand{\cH}{\mathcal H}
\newcommand{\cI}{\mathcal I}
\newcommand{\cJ}{\mathcal J}
\newcommand{\cK}{\mathcal K}
\newcommand{\cL}{\mathcal L}
\newcommand{\cM}{\mathcal M}

\newcommand{\cO}{\mathcal O}
\newcommand{\cP}{\mathcal P}
\newcommand{\cQ}{\mathcal Q}
\newcommand{\cR}{\mathcal R}
\newcommand{\cS}{\mathcal S}

\newcommand{\cU}{\mathcal U}
\newcommand{\cV}{\mathcal V}
\newcommand{\cW}{\mathcal W}

\newcommand{\cZ}{\mathcal Z}
\newcommand{\inj}{\hookrightarrow}
\newcommand{\To}{\longrightarrow}

\newcommand{\emb}{\hookrightarrow}

\renewcommand{\Tilde}{\widetilde}

\renewcommand{\Bar}{\overline}

\newcommand{\cross}{\times}

\newcommand{\imic}{\cong}

\newcommand{\Pic}{\mathop{\mathrm {Pic}}\nolimits}

\newcommand{\orb}{\mathop{\mathrm {orb}}\nolimits}
\newcommand{\rk}{\mathop{\mathrm {rk}}\nolimits}

\newcommand{\bde}{\mathbf e}

\newcommand{\sfa}{\mathsf a}
\newcommand{\sfb}{\mathsf b}
\newcommand{\sfc}{\mathsf c}
\newcommand{\sfd}{\mathsf d}
\newcommand{\sfe}{\mathsf e}
\newcommand{\sff}{\mathsf f}
\newcommand{\sfg}{\mathsf g}
\newcommand{\sfh}{\mathsf h}
\newcommand{\sfi}{\mathsf i}
\newcommand{\sfj}{\mathsf j}
\newcommand{\sfk}{\mathsf k}
\newcommand{\sfl}{\mathsf l}

\newcommand{\sfp}{\mathsf p}
\newcommand{\sfq}{\mathsf q}
\newcommand{\sfr}{\mathsf r}
\newcommand{\sfs}{\mathsf s}
\newcommand{\sft}{\mathsf t}
\newcommand{\sfu}{\mathsf u}
\newcommand{\sfv}{\mathsf v}
\newcommand{\sfw}{\mathsf w}
\newcommand{\sfx}{\mathsf x}
\newcommand{\sfy}{\mathsf y}
\newcommand{\sfz}{\mathsf z}
\newcommand{\qedsymbol}{\mbox{$\Box$}}
\newcommand{\qed}{\unskip\nobreak\hfil\penalty50\hskip1em\hbox{}\nobreak
\hfill\qedsymbol\parfillskip=0pt\finalhyphendemerits=0}
\newenvironment{proof}{\begin{ProofwCaption}{Proof}}{\end{ProofwCaption}}
\newenvironment{ProofwCaption}[1]
 {\addvspace\theorempreskipamount \noindent{\it #1.}\rm}
 {\qed \par \addvspace\theorempostskipamount}

\begin{document}

\title{Numerical obstructions to abelian surfaces in toric Fano 4-folds}
\author{G.K. Sankaran}
\maketitle

On the list of smooth toric Fano 4-folds (\cite{Bat, Sat1}), of which
there are 124, there remain 19 for which it is not known whether there
are embedded abelian surfaces. In order to avoid trivialities, we
consider only totally nondegenerate embedded abelian surfaces, in the
sense of~\cite{San, Sat1}.
\begin{introdefinition}
An embedding $\phi\colon Y \inj X$ from a smooth variety $Y$ to a
smooth toric 4-fold $X$, or its image $\phi(Y)\subset X$, is totally
nondegenerate if $\phi(Y)\cap D\neq \emptyset$ for
any torus-invariant prime divisor~$D$. 
\end{introdefinition}
Results from \cite{HM, Hul, Kaj1, Kaj2, Lan, San, Sat2} about totally
nondegenerate embeddings of abelian surfaces in toric Fano $4$-folds
may be summarised as follows.

\begin{introthm}
  Suppose $X$ is a smooth toric Fano 4-fold. Then $X$ admits a totally
  nondegenerate abelian surface if $X=\PP^4$, if
  $X=\PP^1\cross \PP^3$ (type $\cB_4$), or if $X$ is a product of two
  smooth toric Del Pezzo surfaces (i.e. of type $\cC_4$, $\cD_{13}$,
  $\cH_8$, $\cL_7$, $\cL_8$, $\cL_9$, $\cQ_{10}$, $\cQ_{11}$, $\cK_4$,
  $\cU_5$, $\cS_2\cross \cS_2$, $\cS_2\cross \cS_3$ or $\cS_3\cross
  \cS_3$). Otherwise there is no such embedding, unless possibly $X$
  is of type $\cC_3$, $\cD_7$, $\cD_{10}$, $\cD_{11}$, $\cD_{14}$,
  $\cD_{17}$, $\cD_{18}$, $\cG_3$, $\cG_4$, $\cG_5$, $\cL_{11}$,
  $\cL_{13}$, $\cI_9$, $\cQ_{16}$, $\cU_8$, $\cV_4$, $\cW$, $\cZ_1$ or
  $\cZ_2$.
\end{introthm}

The notation for the types is taken from \cite{Sat1}, which apart from
the addition of the missing case $\cW$ is the same as in \cite{Bat}. 

It is claimed in \cite{San} that there are no totally nondegenerately
embedded abelian surfaces in toric Fano 4-folds of types $\cC_1$ or
$\cC_2$ but that there are in the case of type $\cC_3$. However, it
was discovered by T.~Kajiwara that there is an error in the
calculations in \cite{San}: in the case of $\cC_3$ it is claimed,
wrongly, that a certain class satisfies the condition given by the
double-point formula (Lemma~\ref{double} in this paper). I take this
opportunity to thank Professor Kajiwara for pointing this out to me. 

In view of this quantity of previous work it is perhaps surprising
that more can be achieved by completely elementary methods.
Nevertheless, in this paper we eliminate or restrict the remaining
cases by examining the numerical conditions imposed on a class
$\alpha$ in $A^2(X)$ by the condition that it should contain a minimal
abelian surface. The methods are entirely elementary but they do rely
on some small computer calculations.

There are also some geometric conditions. A simple example is that if
$f\colon X \to S$ is a map to a rational surface then $f\circ
\phi\colon A\to S$ cannot have degree~$1$. We make some remarks on
such conditions but we do not attempt to exploit them systematically here. 

The most conclusive results of this paper are as follows. 

\begin{introthm}
  Suppose $X$ is a smooth toric Fano 4-fold. Then $X$ does not contain
  a totally nondegenerate abelian surface if it is of type
  $\cC_1$, $\cC_2$, $\cC_3$, $\cD_{10}$ or $\cD_{18}$. 
\end{introthm}

In many other cases we can prove that if a class $\alpha\in A^2(X)$
contains a totally nodegenerate abelian surface then it satisfies
certain restrictive conditions.
\begin{itemize}
\item If $X$ is of type $\cD_7$, $\cG_4$, $\cI_9$, $\cL_{13}$,
  $\cQ_{16}$, $\cU_8$ or $\cZ_1$, then the class $\alpha$ of such a
  totally nondegenerate abelian surface in $X$ must belong to a short
  finite list.
\item If $X$ is of class $\cD_{14}$ or $\cD_{17}$ then $\alpha$ belongs to a
  short finite list or satisfies other strong conditions. 
\item If $X$ is of type $\cD_{11}$, $\cG_3$, $\cG_5$, $\cW$ or $\cZ_2$ then the
  coefficients of the class $\alpha$ satisfy a condition modulo~$3$
  and some other weak conditions. 
\item If $X$ is of type $\cL_{11}$ or $\cV_4$ then we obtain only weak
  conditions on the class~$\alpha$. 
\end{itemize}
The precise statements may be found in
sections~\ref{typeC}--\ref{typeother}. 

\section{Preliminaries}\label{tools}

In this section, unless otherwise stated, $X$ is any smooth projective
toric $4$-fold. 

Suppose that $A$ is a (minimal) abelian surface, and that there is a
totally nondegenerate embedding $\phi\colon A \to X$.  Then $\phi$
determines a class $\alpha=[\phi(A)]\in A^2(X)\imic H^4(X,\ZZ)$. There
are strong restrictions on $\alpha$. If $X$ is one of the 124 smooth
toric Fano 4-folds listed in~\cite{Sat1}, then in most cases the
restrictions are so strong that no such $\phi$ can exist. 

Let $D_1,\ldots,D_r$ be the torus-invariant prime divisors on $X$, and
let $\Gamma_i$ be the corresponding curve on $A$, so
$\cO_A(\Gamma_i)=\phi^*\big(\cO_X(D_i)\big)$. We denote by $L$ the
quadratic form given by $L_{ij}=(\Gamma_i\cdot\Gamma_j)_A$: note that
$L_{ij}=D_iD_j\alpha\in A^4(X)\imic \ZZ$. Since $X$ is assumed to be
smooth and projective, it follows easily (see~\cite[2.1]{Od}) that
there exist $\lambda_i\ge 0$ such that $\sum \lambda_iD_i$ is ample.

\begin{lemma}
The quadratic form $L$ has the following properties. 
\begin{itemize}
\item[(i)] The rank of $L$ is $e\le 4$ and the signature is $(1,e-1)$. 
\item[(ii)] $L$ is an even form (in particular, $L_{ii}$ is even), and
  $L_{ij}\ge 0$ for all $i$, $j$. 
\item[(iii)] If $L_{ij}=0$ then $L_{ii}=L_{jj}=0$. 
\item[(iv)] For any $i$, there exists $j$ such that $L_{ij}>0$. 
\item[(v)] If $L_{ij}=0$ then $(L_{i1}\colon\ldots\colon L_{ir})=
  (L_{j1}\colon\ldots\colon L_{jr})\in\PP^{r-1}_\QQ$. In particular
  $L_{ik}L_{jl}=L_{il}L_{jk}$ for all $k$ and $l$. 
\item[(vi)]\label{degree1}   If $L_{ii}=0$ and $L_{ij}=1$ then $2L_{jk}\ge L_{jj}L_{ik}$ for
  all~$k$. 
\end{itemize}
\end{lemma}
\begin{proof}
  \begin{itemize}
  \item[(i)] 
  $L$ is isomorphic to a sublattice of $H^2(A,\ZZ)$ and is generated
  by algebraic cycles. But $\rho=\rk\Pic A\le 4$ (for any abelian
  surface $A$), see \cite{LB}, and the signature of the intersection
  form on $\Pic A$ is $(1,\rho-1)$ by the Hodge index theorem. 
  Moreover, $\sum \lambda_iD_i$ is ample on $X$ and hence $L$ is not
  negative semidefinite. 
\item[(ii)]
  These are standard properties of the intersection form on an abelian
  surface. 
\item[(iii)]
  See for instance \cite{Sat2}, or note that by the Nakai criterion
  \cite[Corollary 4.3.3]{LB} an effective curve with
  positive square on an abelian surface is ample. 
\item[(iv)]
  Since $\sum \lambda_i D_i$ is ample on $X$ and $\phi$ is an embedding, $\sum
  \lambda_i\Gamma_i$ is ample on $A$. Therefore $0<\Gamma_i\cdot (\sum
  \lambda_j \Gamma_j)=\sum_j \lambda_jL_{ij}$ so some $L_{ij}$ is positive. 
\item[(v)]
  Since $L_{ij}=0$ we have by Lemma~1.1(iii) that
  $\Gamma_i^2=\Gamma_j^2=0$. Let $\Gamma$ be a reduced irreducible
  component of the effective (possibly nonreduced) curve $\Gamma_i$. 
  If $\Gamma'$ is any reduced irreducible component of $\Gamma_i$ (or
  $\Gamma_j$) we have $0\le \Gamma\cdot\Gamma' \le \Gamma\cdot
  \Gamma_i \le \Gamma_i^2=0$, so in particular $\Gamma^2=0$. Since $A$
  contains no rational curves, $\Gamma$ must be an elliptic curve and
  we consider the exact sequence
$$
0\To \Gamma \To A \To E \To 0. 
$$
Every component $\Gamma'$ as above is contained in a fibre of $A\to
E$, since $\Gamma'\cdot \Gamma=0$, so the numerical classes
$[\Gamma_i]$ and $[\Gamma_j]$ are both integer multiples of $[\Gamma]$
and therefore rational multiples of each other. 
\item[(vi)]
  As in the proof of Lemma~1.1(v),
  $[\Gamma_i]=n[\Gamma]$ for some elliptic curve $\Gamma$ and
  some $n\in \NN$. But $n|L_{ij}$ so $n=1$ and
  $\Gamma_i$ is a reduced elliptic curve. So we have an exact sequence
$$
0\To \Gamma_i \To A \To E \To 0
$$ 
where $E$ is an elliptic curve: in fact $A\imic \Gamma_i\times E$. But
$\Gamma_j\to E$ is generically finite of degree~$1$, so
$[\Gamma_j]=[E+m\Gamma_i]$ for some $m\in \ZZ$ (identifying $E$
with the zero section). Since $E\cdot\Gamma_i=1$ and
$E^2=\Gamma_i^2=0$, we have $L_{jj}=\Gamma_j^2=2m$, and hence
$$
0\le 2\Gamma_k\cdot E=2L_{jk}-2m L_{ik}=2L_{jk}-L_{jj}L_{ik}. 
$$
as claimed. 
\end{itemize}
\end{proof}

In particular, if $L_{ii}=0$ then $\Gamma_i$ is a union of disjoint
smooth genus~$1$ curves (with multiplicity), all translates of one
another. 

Another constraint comes from the double-point formula. 

\begin{lemma}\label{double}
  If $\alpha\in A^2(X)$ is the class of a smooth minimal abelian
  surface embedded in $X$ by $\phi\colon A\emb X$ then
  $\alpha^2-\alpha\cdot c_2(X)=0$. 
\end{lemma}

\begin{proof}
  The double-point scheme is given (\cite[Theorem 9.3]{Ful}) by
\begin{eqnarray*}
  \DD(\phi)&=&\phi^*\phi_*[A]-\big(c(\phi^*T_X)c(T_A)^{-1}\big)_2\cap[A]\\
  &=&\phi^*\phi_*[A]-c_2(\phi^*T_X)\cap[A]\\
  &=&\big([\phi(A)]-c_2(T_X)\big)\cap[\phi(A)]. 
\end{eqnarray*}
Since $\phi(A)$ is smooth $\DD(\phi)$ has length zero, i.e.{}
$\big([\phi(A)]-c_2(T_X)\big)\cdot[\phi(A)]=0$ in $A^4(X)\imic\ZZ$,
which is what is claimed. 
\end{proof}

If there is a toric map from $X$ to a curve or a surface
it is sometimes possible to use it to obtain further restrictions
on~$L_{ij}$. This happens for a few of the smooth toric Fano $4$-folds. 

\section{Notation and methods}\label{general}

For the rest of the paper, we let $X=T_N\emb(\Sigma)$ be a smooth
toric Fano 4-fold given by a fan $\Sigma$. Following \cite{Bat} we put
$\Sigma^{(1)}=\{\tau_i\in \Sigma\mid \dim \tau_i=1\}$ and $G(\Sigma) =
\{x_i\mid \tau_i \in \Sigma^{(1)}\}$, where $x_i$ is the unique
element of $\tau_i$ that generates the lattice $\tau_i\cap N$. The fan
$\Sigma$ is determined up to an integral change of basis by the
primitive relations (see \cite{Bat}), which are tabulated in
\cite[1.9]{Sat3}. We adopt the convention that
$D_i=\Bar{\orb(\tau_i)}$ denotes the divisor corresponding to $x_i\in
G(\Sigma)$ as listed in \cite{Sat3}. 

For each type of Fano 4-fold we use {\it Macaulay2\/} to calculate the
Stanley-Reisner ring, i.e. $H^{2*}(X,\ZZ)=A^*(X)$ with the
intersection product. It is generated in degree~$1$, by
$D_1,\ldots,D_r$. Note that $\rk\Pic(X)=r-4$. There are linear
relations $\sum_i m\big(\bde(\tau_i)\big)D_i=0$ for $m\in N^\vee$ (of
course it is enough to take $m$ in a $\ZZ$-basis of $N^\vee$) and
multiplicative relations $\prod_{i\in I}D_i=0$ if $\sum_{i\in I}\tau_i
\not\in \Sigma$. All these relations are easily computed from the
primitive relations. 

We also calculate a $\ZZ$-basis $\beta_1,\ldots,\beta_s$ for $A^2(X)$,
in which each $\beta_k$ is of the form $\left[D_iD_j\right]$ for some
$1\le i\le j\le r$. This is simply a basis for the degree~$2$ part of
the Stanley-Reisner ring. For a class $\alpha=\sum a_k\beta_k\in
A^2(X)$ we calculate the intersection matrix $L_{ij}=D_iD_j\alpha\in
A^4(X)\imic\ZZ$. 

For ease of reading we write $\sfa$, $\sfb$, $\sfc$, etc.{} instead of
$a_1$, $a_2$, $a_3$ etc. 

We do not always need the whole of $L$. If $i<j$ and $D_i\equiv D_j$
(numerical equivalence) then $L_{ik}=L_{jk}$ for all $k$ and we may
omit the redundant $i$th row and column altogether. We denote by
$\Lambda$ the submatrix of $L$ thus obtained. Thus $\Lambda$ is the
restriction of $L$ to the span of all $D_i$ that are not numerically
equivalent to any $D_j$ with $j>i$, and $\Lambda_{\mu\nu}=L_{ij}$ if
$D_i$ and $D_j$ are in the $\mu$th and $\nu$th (in the revlex order)
numerical equivalence class respectively. 

To apply this to particular cases we assume that $\alpha$ is the class
of a smooth minimal abelian subvariety of $X$, and deduce conditions
on $\sfa, \sfb,\dots$ which in many cases lead to a
contradiction. There is some further computation involved, which was
done using {\it Maple}. 

\section{Types $\cC_1$, $\cC_2$ and $\cC_3$}\label{typeC}

In this section we re-examine the cases $\cC_1$, $\cC_2$ and $\cC_3$. 
In particular, we correct the error in~\cite{San} described in the
introduction above, by giving a new analysis of the case $\cC_3$. For
$\cC_1$ and $\cC_2$ the results are as in \cite{San}, but we reprove
them here more concisely and in the notation used in the rest of this
paper. 

In these cases $\rk\Pic(X)=2$. The fans of the Fano $4$-folds of these
types are given by the primitive relations shown in the table. 
\begin{center}
\begin{tabular}
{|c||c|c|c||}
\hline
&$\cC_1$&$\cC_2$&$\cC_3$\\
\hline
$x_1+x_2+x_3=$&$0$&$0$&$0$\\
\hline
$x_4+x_5+x_6=$&$2x_1$&$x_1$&$x_1+x_2$\\
\hline
\end{tabular}
\end{center}
For $X$ of type $\cC_1$ or $\cC_2$ we have $D_2\equiv D_3$ and
$D_4\equiv D_5$, and a basis for $A^2(X)$ is $\beta_1=D_3^2$,
$\beta_2=D_3D_6$, $\beta_3=D_6^2$. 

\begin{proposition}\label{C1C2} 
There are no totally nondegenerate abelian surfaces in toric Fano
4-folds of types $\cC_1$ or $\cC_2$. 
\end{proposition}

\begin{proof}
  We need the observation from \cite{San} that (in our present
  notation) $\sfa\ge 6$ (and is even). In fact
  $\sfa=\Lambda_{33}=D_6^2\alpha$ in both cases, and $D_4$, $D_5$ and
  $D_6$ are the pull-backs of lines in $\PP^2$ under a projection
  $p\colon X\to \PP^2$ . These Fano 4-folds are the projectivisations
  of toric (hence decomposible) rank~$3$ vector bundles on $\PP^2$,
  namely $\PP\big(\cO_{\PP^2} \oplus \cO_{\PP^2} \oplus
  \cO_{\PP^2}(2)\big)$ for $\cC_1$ and $\PP\big(\cO_{\PP^2} \oplus
  \cO_{\PP^2} \oplus \cO_{\PP^2}(1)\big)$ for $\cC_2$, and $p$ is the
  projection. Hence if $A\subset X$ is an abelian surface with
  $[A]=\alpha$, and $p|_A\colon A\to \PP^2$ is surjective, then by
  Riemann-Roch
$$
3=h^0\big(\cO_{\PP^2}(1)\big)\le h^0\big(p|_A^*\cO_{\PP^2}(1)\big)=
\frac{1}{2}D_6^2\alpha=\frac{1}{2}\sfa. 
$$
On the other hand, if $h^0\big(\cO_A(D_6)\big)\le 2$ then $p|_A$ is
given by a subsystem of the linear system $|D_6|$ and hence it maps
$A$ onto a line in $\PP^2$. But the inverse image of this line is a
smooth toric 3-fold and hence contains no abelian surface (see for
instance \cite{Kaj2}). 
\medskip

For $\cC_1$ we find that $\Lambda_{12}=\sfc$ and $\Lambda_{13}=\sfb$, so $\sfb\ge
0$ and $\sfc \ge 0$ by Lemma~1.1(ii).  But the double-point
formula (Lemma~\ref{double}) gives
\begin{eqnarray*}
0&=&  4\sfa^2+4\sfa\sfb+\sfb^2+2\sfa\sfc-19\sfa-11\sfb-3\sfc\\
&=& (\sfa+\sfb)^2+(2\sfa-11)\sfb+(\sfa-3)\sfc+(3\sfa-19+\sfc)\sfa
\end{eqnarray*}
which is obviously positive if $\sfa\ge 6$. 
\medskip

For $\cC_2$ we again have $0\le \Lambda_{12}=\sfc$ and
$0\le \Lambda_{13}=\sfb$. The double-point formula gives
\begin{eqnarray*}
0&=&\sfa^2+2\sfa\sfb+\sfb^2+2\sfa\sfc-10\sfa-10\sfb-3\sfc\\
&=&(\sfa-10)\sfa+\sfb^2+(2\sfa-10)\sfb+(2\sfa-3)\sfc
\end{eqnarray*}
which is positive if $\sfa\ge 10$. But $0\le \Lambda_{11}=\sfc-\sfb$, so we get
$$
0=-16+6\sfb+\sfb^2+13\sfc\ge\sfb^2+19\sfb-16
$$ 
for $\sfa=8$, and for $\sfa=6$ 
$$
0=-24+2\sfb+\sfb^2+9\sfc\ge
\sfb^2+11\sfb-24. 
$$
Neither of these has a solution with $\sfb$, $\sfc$
non-negative integers. 
\end{proof}

For $X$ of type $\cC_3$ we have $D_1\equiv D_2$ and $D_4\equiv
D_5\equiv D_6$. A basis for $A^2(X)$ is given by $\beta_1=D_3^2$,
$\beta_2= D_3D_6$ and $\beta_3= D_6^2$. 

\begin{proposition}\label{C3} 
There are no totally nondegenerate abelian surfaces in a toric Fano
4-fold of type $\cC_3$. 
\end{proposition}

\begin{proof} The computation gives
\begin{equation*}
\Lambda=\begin{pmatrix}
 \sfc&\sfx+\sfc&\sfx\\
\sfx+\sfc&\sfa+2\sfx+\sfc&\sfa+\sfx\\
\sfx&\sfa+\sfx&\sfa
\end{pmatrix}
\end{equation*}
where $\sfx=\sfa+\sfb$. As $\Lambda$ has rank at most~$2$, the Hodge
index theorem gives $\sfx^2\ge \sfa\sfc$, and $\sfx>0$ since otherwise
a row of $\Lambda$ vanishes. We may assume $\sfa\ge 6$, as in
Proposition~\ref{C1C2}. 

The double-point formula gives
\begin{eqnarray}\label{dbpC3}
0&=&3\sfa^2+4\sfa\sfb+\sfb^2+2\sfa\sfc-17\sfa-11\sfb-3\sfc\nonumber\\
&=&2\sfa\sfx+\sfx^2+2\sfa\sfc -6\sfa-11\sfx-3\sfc\\
&\ge&\sfx+3\sfa\sfc-6\sfa-3\sfc\nonumber\\  
&>&3((\sfa-1)(\sfc-2)-2).\nonumber 
\end{eqnarray}
Hence $\sfc=0$ or $\sfc=2$. Since $\sfx\le 11$ immediately from
equation~(\ref{dbpC3}), it is simple to check that the only
possibility is $\sfc=0$, $\sfx=4$, $\sfa=14$. 

Suppose this occurs. Then on $\phi(A)$ we have $D_1^2=0$, so $D_1$ is
a union of elliptic curves in $A$. It cannot be a single reduced
elliptic curve because it moves in the linear system $D_1+tD_2$, so
$A$ has a $(1,7)$-polarisation given by $D_6$, though which $p\colon
A\to \PP^2$ factors, and an elliptic curve $C$ on which $D_6$ has
degree $1$ or $2$. The former case is evidently impossible since then
we do not even get a morphism to $\PP^2$ defined on $C$. 

In the latter case, the linear system $\vert D_6\vert$ is not very
ample: $(A, \cO_A(D_6))$ is bielliptic and $\phi_{\vert D_6\vert}$
maps $A$ 2-to-1 onto its image~$S$, which is an elliptic ruled surface
in $\PP^6$. The fibres of this ruled surface are the Kummer curves of
the translates of $C$. Hence pairs of points in each such translate
are identified both by the subsystem $\vert p^*(\cO_{\PP^2}(1))\vert
\subset \vert D_6\vert$ and by the linear system $\vert D_1\vert$. In
particular the induced map $A\to \PP^2\cross\PP^1$ factors
through~$S$. According to \cite[Section 2]{San}, to specify a map
to~$\cC_3$ we have to choose also a subsystem of $\vert D_1+D_6\vert$
given by $D_3$, i.e. a section $\sigma\in
H^0\big(\cO_A(D_1+D_6)\big)$. 

The curve $\Gamma_3=D_3\vert_A=\{\sigma=0\}\subset A$ has genus given by
$2p_a(\Gamma_3)-2= (D_1+D_6)^2=22$ so $p_a(\Gamma_3)=12$. Consider the
image $\Delta_3=\phi_{\vert D_6\vert}(\Gamma_3)\subset S$. It has
$\Delta_3^2=\Gamma_3^2=22$ and $\Delta_3.K_S\neq 0$ since $-K_S$ is
effective and $\Delta_3$ moves. So $p_a(\Delta_3)\neq p_a(\Gamma_3)$
so $\Delta_3$ is not isomorphic to $\Gamma_3$; hence there are two
(possibly infinitely near) points of $\Gamma_3$ that are identified by
$A\to\PP^2\cross\PP^1$. But they are also not separated by $\sigma$,
so the map $A\to \cC_3$ is not an isomorphism on~$\Gamma_3$. 
\end{proof}

\section{Types $\cD_7$, $\cD_{10}$, $\cD_{11}$, $\cD_{14}$, 
$\cD_{17}$ and $\cD_{18}$}\label{typeD}

In these cases $\rk\Pic(X)=3$. We obtain strong restrictions on
$\alpha$ except in case $\cD_{11}$. The fans of the Fano $4$-folds of
these types are given by the primitive relations shown in the table. 
\begin{center}
\begin{tabular}
{|c||c|c|c|c|c|c||}
\hline
&$\cD_7$&$\cD_{10}$&$\cD_{11}$&$\cD_{14}$&$\cD_{17}$&$\cD_{18}$\\
\hline
$x_1+x_2+x_3=$&$0$&$x_6$&$0$&$0$&$0$&$2x_7$\\
\hline
$x_4+x_5=$&$x_1$&$x_1$&$x_1$&$0$&$x_1$&$x_6$\\
\hline
$x_6+x_7=$&$x_1$&$0$&$x_4$&$x_1$&$x_2$&$0$\\
\hline
\end{tabular}
\end{center}
For $X$ of type $\cD_7$ we have $D_2\equiv D_3$, $D_4\equiv D_5$ and
$D_6\equiv D_7$. We choose the basis $\beta_1=D_3^2$, $\beta_2=
D_3D_5$, $\beta_3= D_3D_7$, $\beta_4= D_5D_7$ for $A^2(X)$. 
Interchanging $D_5$ and $D_7$ if necessary, we may assume that
$\sfb\ge \sfc$. 
\begin{proposition}\label{D7}
  The class of any totally nondegenerate abelian surface in a Fano
  $4$-fold of type $\cD_7$ satisfies one of
  $(\sfa,\sfb,\sfc,\sfd)=(3,4,1,7)$, $(3,12,0,12)$, $(4,4,0,4)$ or
  $(5,2,0,2)$. 
\end{proposition}

\begin{proof}
  The $4$-fold $X$ of type $\cD_7$ is a $\PP^2$-bundle over
  $\PP=\PP^1\times \PP^1$. If $p\colon X\to \PP$ is the projection
  then $D_5$ and $D_7$ are the pullbacks of lines in the two rulings. 
  The computation gives $\Lambda_{34}=L_{57}=\sfa$, so if $\sfa\neq 0$ then
$$
3=h^0\big(\cO_\PP(1,1)\big)\le h^0\big(\cO_A(D_5+D_7)\big)=\sfa
$$
by Riemann-Roch. The double-point formula gives 
\begin{equation}\label{dbpD7}
(2\sfa^2-14\sfa)+(2\sfa-7)(\sfb+\sfc)+(2\sfa-3)\sfd=0
\end{equation}
so evidently $0<\sfa\le 6$. Moreover, $\Lambda_{11}=\sfd-\sfb-\sfc$
and $\Lambda_{13}=\sfc$, so by Lemma~1.1(ii) $\sfd\ge
\sfb+\sfc$ and $\sfc\ge 0$, and if $\sfc=0$ then $\sfb=\sfd$ by
Lemma~1.1(iii). 

With these restrictions it is easy to see that the solutions to
equation~(\ref{dbpD7}) are $(\sfa,\sfb,\sfc,\sfd)=(3,1,1,8)$,
$(3,4,1,7)$, $(3,12,0,12)$, $(4,1,1,4)$, $(4,4,0,4)$ and $(5,2,0,2)$. 
However the cases $(3,1,1,8)$ and $(4,1,1,4)$ are excluded by
Lemma~1.1(vi). 
\end{proof}

For $X$ of type $\cD_{10}$ we have $D_2\equiv D_3$ and $D_4\equiv
D_5$. We choose the basis $\beta_1=D_5D_6$, $\beta_2=D_5D_7$,
$\beta_3=D_6^2$, $\beta_4=D_7^2$ for $A^2(X)$. 

\begin{proposition}\label{D10}
  There are no totally nondegenerate abelian surfaces in a toric Fano
  $4$-fold of type $\cD_{10}$. 
\end{proposition}

\begin{proof}
The computation yields
\begin{equation*}
\Lambda_{4*}= \begin{pmatrix}
      -\sfa&
      -\sfa+\sfc&
      \sfc&
      \sfa-\sfc&
      0\\
      \end{pmatrix}. 
\end{equation*}
So $-\sfa+\sfc\ge 0$ and $\sfa-\sfc\ge 0$ by Lemma~1.1(ii),
so $\sfa=\sfc$; but again by Lemma~1.1(ii), $-\sfa\ge 0$ and
$\sfc\ge 0$ so $\sfa=\sfc=0$. But now $\Lambda_{4j}=0$ for all~$j$, so
by Lemma~1.1(iv) no totally nondegenerate embedding exists. 
\end{proof}

For $\cD_{11}$ the numerical conditions are rather
weak. We have $D_2\equiv D_3$ and $D_6=D_7$, and we choose the basis
$\beta_1=D_3^2$, $\beta_2= D_3D_5$, $\beta_3=D_3D_7$,
$\beta_4=D_5D_7$. 
\begin{proposition}\label{D11}
  The class of any totally nondegenerate abelian surface in a toric Fano
  $4$-fold of type $\cD_{11}$ satisfies $\sfa=0$, $\sfd-\sfb-\sfc\ge 0$ and
  $3\sfd=\sfb^2+2\sfb\sfc-10\sfb-7\sfc$. 
\end{proposition}

\begin{proof}
  The computation gives $\Lambda_{55}=-\Lambda_{33}=\sfa$, so
  $\sfa=0$, and $\Lambda_{11}=\sfd-\sfb-\sfc$, and the condition
  $3\sfd=\sfb^2+2\sfb\sfc-10\sfb-7\sfc$ comes from the double-point
  formula. 
\end{proof}

For $\cD_{14}$ the numerical conditions are fairly restrictive and the
existence of a map to $\PP^1\times\PP^1$ reduces the possibilities
still further. We have $D_2\equiv D_3$, $D_4\equiv D_5$ and $D_6\equiv
D_7$ and we choose the basis $\beta_1=D_3^2$, $\beta_2= D_3D_5$,
$\beta_3= D_3D_7$, $\beta_4= D_5D_7$. 

\begin{proposition}\label{D14}
  The class of any totally nondegenerate abelian surface in a toric
  Fano $4$-fold of type $\cD_{14}$ satisfies one of $\sfa\le 3$,
  $\sfb\le 2$ or $(\sfa,\sfb,\sfc,\sfd)=(4,4,4,4)$ or $(5,4,0,4)$. 
\end{proposition}
\begin{proof}
The computation yields $\Lambda_{11}=\sfd-\sfb$, $\Lambda_{13}=\sfc$,
$\Lambda_{14}=\sfb$ and $\Lambda_{34}=\sfa$, and the double-point formula
gives
\begin{eqnarray*}
0&=&2\sfa\sfb+2\sfb\sfc+2\sfa\sfd-8\sfa-7\sfb-6\sfc-3\sfd\\
&=&4\sfa\sfb+2\sfb\sfc+2\sfa(\sfd-\sfb)-8\sfa-10\sfb-6\sfc-3(\sfd-\sfb)\\
&=&\sfa(2\sfb-8)+\sfb(2\sfa-10)+\sfc(2\sfb-6)+\Lambda_{11}(2\sfa-3). 
\end{eqnarray*}
If we assume that $\sfa\ge 4$ and $\sfb\ge 3$, this immediately
gives $\sfb\le 4$ or $\sfa\le 5$. Moreover the Hodge index theorem
gives $2\sfb\sfc+\sfa(\sfd-\sfb)\le 0$. It is easy to check that the
only solutions to this are those claimed. 
\end{proof}

For $\cD_{17}$ the position is similar to that for $\cD_{14}$. For $X$
of type $\cD_{17}$ we have $D_4\equiv D_5$ and $D_6\equiv D_7$ and we
choose the basis $\beta_1=D_3^2$, $\beta_2= D_3D_5$, $\beta_3=
D_3D_7$, $\beta_4= D_5D_7$. Interchanging $D_5$ and $D_7$ if
necessary, we may assume that $\sfb\ge \sfc$, and we write
$\sfx=\sfb-\sfc\ge 0$. 

\begin{proposition}\label{D17}
  The class of any totally nondegenerate abelian surface in a toric
  Fano $4$-fold of type $\cD_{17}$ satisfies one of $\sfa=0$ or
  $(\sfa,\sfb,\sfc,\sfd)=(3,3,3,5)$, $(3,4,1,8)$, $(3,9,0,12)$,
  $(4,2,0,6)$, $(4,2,1,5)$ or $(4,2,2,4)$. 
\end{proposition}
\begin{proof}
Again $\sfa=\Lambda_{45}$ is the degree of a map $A\to \PP^1\times \PP^1$
and hence $\sfa=0$ or $\sfa\ge 3$. We also have $\Lambda_{12}=\sfd$,
$\Lambda_{15}=\sfb$ and $\Lambda_{24}=\sfc$. Moreover
$\Lambda_{11}+\Lambda_{22}=2(\sfd-\sfa)$ so $\sfd\ge\sfa$, and we
write $\sfy=\sfd-\sfa$. Now the double-point formula gives
\begin{eqnarray*}
  0&=&\sfa^2+2\sfa\sfb+2\sfa\sfc+2\sfb\sfc+2\sfa\sfd-12\sfa-7\sfb-7\sfc-3\sfd\\
&=&(3\sfa-15)\sfa+(4\sfa-14)\sfc+(2\sfa-7)\sfx+(2\sfa-3)\sfy+2\sfc^2+2\sfc\sfx
\end{eqnarray*}
which immediately implies that $\sfa\le 4$. Together with the Hodge
index theorem, which gives
$\sfa^2+\sfa\sfc-\sfa\sfd+\sfa\sfb+2\sfb\sfc\ge 0$ for $\sfa\neq 0$,
this is enough to yield the result by a simple calculation. 
\end{proof}

Finally, the case $\cD_{18}$ may also be excluded altogether. In this
case we have $D_1\equiv D_2\equiv D_3$ and $D_4\equiv D_5$, and we choose
the basis $\beta_1=D_3D_5$, $\beta_2=D_3D_7$, $\beta_3= D_5D_7$,
$\beta_4=D_3^2$. 

\begin{proposition}\label{D18}
  There are no totally nondegenerate abelian surfaces in a toric Fano
  $4$-fold of type $\cD_{18}$. 
\end{proposition}

\begin{proof}
We find
\begin{equation*}
\Lambda= \begin{pmatrix}
      \sfc&
      \sfb&
      \sfa&
      \sfa+\sfb-2 \sfc\\
      \sfb&
      0&
      \sfd&
      -2 \sfb+\sfd\\
      \sfa&
      \sfd&
      2 \sfa-\sfd&
      0\\
      \sfa+\sfb-2 \sfc&
      -2 \sfb+\sfd&
      0&
      -2 \sfa-4 \sfb+4 \sfc+\sfd\\
      \end{pmatrix}, 
\end{equation*}
so by Lemma~1.1(iii) we have $\sfd=2\sfa$ and $\sfb=\sfc$. 
The double-point formula gives
\begin{eqnarray}\label{dbpD18}
0&=&
2\sfa\sfb+\sfb^2-4\sfb\sfc+2\sfc\sfd-6\sfa-5\sfb+3\sfc-4\sfd\nonumber\\
&=&6\sfa \sfb -3\sfb^2-14\sfa-2\sfb
\end{eqnarray}
and by Lemma~1.1(ii) $\sfa\ge 0$, $\sfb\ge 0$ and $0\le
\sfa+\sfb-2\sfc=\sfa-\sfb$. If $\sfb=0$ then $\sfa=0$ by
equation~(\ref{dbpD18}), which is impossible by Lemma~1.1(iv)
because $\Lambda_{3*}=0$, and if $\sfa=\sfb$ then $\Lambda_{14}=0$ so $\Lambda_{11}=0$ by
Lemma~1.1(iii); but then $\sfb=0$. Hence $0<\sfb<\sfa$, and $\sfb=\Lambda_{11}$ is
even by Lemma~1.1(ii). So
$$
2\le \sfb <\sfa = (3\sfb^2+2\sfb)/(6\sfb-14)
$$
from equation~(\ref{dbpD18}), since $6\sfb-14\neq 0$. So either $6\sfb-14< 0$,
i.e.{} $\sfb=2$, or $\sfb(6\sfb-14)<3\sfb^2+2\sfb$, which gives
$\sfb\le 16/3 <6$. So $\sfb=2$ or $\sfb=4$; but if $\sfb=2$ then
$\sfa=-8<0$, and if $\sfb=4$ then $\sfa=28/5\not\in \ZZ$. 
\end{proof}

In Proposition~\ref{D14} the conditions $\sfa\ge 4$ and $\sfb\ge 3$
are not arbitrary: they correspond to geometric conditions on~$A$. 
Exactly as in Proposition~\ref{D7}, there is a map to $\PP^1\times
\PP^1$ whose restriction is of degree~$\sfa=L_{57}$. Therefore in any
case $\sfa=0$ or $\sfa\ge 3$. There is also a map from a blow-up
$\tilde X$ of $X$ to the Hirzebruch surface $\FF_1$ which is of degree
$\sfb$ on the proper transform $\tilde A$ of $A$. We obtain the
blow-up by introducing the ray $\RR_+(-x_2)$ into $\Sigma^{(1)}$: the
map to $\FF_1$ is induced by projection to the plane spanned by $x_1$
and $x_7$. Since $\FF_1$ is rational and $\tilde A$ is not, $\sfb\neq 1$. 

However, if $\sfa=0$ or $\sfb=0$ then $A$ or $\tilde A$ maps
to a curve in $\PP^1\times \PP^1$ or $\FF_1$; if $\sfa=3$ then $A$ is
a triple cover of $\PP^1\times \PP^1$ and if $\sfb=2$ then $\tilde A$
is a double cover of $\FF_1$. All these are strong geometric
constraints. For instance, it is easy to see, using the results of
Miranda ~\cite{Mir} on triple covers, that $A\to \PP^1\times \PP^1$
cannot have general branching behaviour. A similar remark applies to
the cases in Proposition~\ref{D7} and Proposition~\ref{D17} for which $\sfa=3$. 

\section{Types $\cG_3$, $\cG_4$  and $\cG_5$}\label{typeG}

In these cases $\rk\Pic X=3$. For $\cG_3$ and $\cG_5$ the numerical
conditions are quite weak, but for $\cG_4$ they are very restrictive. 
The fans of the Fano $4$-folds of these types are given by the
primitive relations shown in the table. 
\begin{center}
\begin{tabular}
{|c||c|c|c||}
\hline
&$\cG_3$&$\cG_4$&$\cG_5$\\
\hline
$x_1+x_7=$&$0$&$x_4$&$x_4$\\
\hline
$x_2+x_3+x_4=$&$0$&$x_7$&$x_7$\\
\hline
$x_4+x_5+x_6=$&$x_1$&$x_1+x_2$&$0$\\
\hline
$x_5+x_6+x_7=$&$x_2+x_3$&$x_2$&$x_2+x_3$\\
\hline
$x_1+x_2+x_3=$&$x_5+x_6$&$0$&$0$\\
\hline
\end{tabular}
\end{center}
For $X$ of type $\cG_3$ we have $D_2\equiv D_3$ and $D_5\equiv D_6$
and we choose the basis $\beta_1=D_4^2$, $\beta_2=D_4D_6$,
$\beta_3=D_4D_7$, $\beta_4=D_6^2$, $\beta_5=D_7^2$. 

\begin{proposition}\label{G3}
   The class of any totally nondegenerate abelian surface in a toric Fano
  $4$-fold of type $\cG_3$ satisfies $\sfa=\sfd=0$, $\sfb>0$, $\sfc>0$ and
  $3\sfe=2\sfb\sfc+\sfc^2-6\sfb-11\sfc\ge 0$. 
\end{proposition}
\begin{proof}
We have $\Lambda_{15}=0$ and $\Lambda_{15}=\sfd-\sfa$,
$\Lambda_{55}=\sfa$. The formula for $\sfe$ in terms of $\sfb$ and
$\sfc$ comes from the double-point formula. Each of $\sfb$,
$\sfc$ and $\sfe$ occurs as  an entry in $\Lambda$, so they are
non-negative by Lemma~1.1(ii) and $\sfb$ and $\sfc$ are the only
non-zero values in some row, and hence positive by Lemma~1.1(iv). 
\end{proof}

In the case of $\cG_4$ we have $D_5\equiv D_6$ and we choose the basis
$\beta_1=D_3^2$, $\beta_2=D_3D_6$, $\beta_3= D_6^2$, $\beta_4=D_6D_7$,
$\beta_5=D_7^2$. 

\begin{proposition}\label{G4}
  If $\sfa\beta_1+\cdots+\sfe\beta_5$ is the class of a totally
  nondegenerate abelian surface in a toric Fano $4$-fold of
  type $\cG_4$ then $(\sfa,\sfb,\sfc, \sfd, \sfe)$ is one of
$$
(2,6,0,8,-1), (12,-8,0,12,-10),
$$
$$
(12,-8,0,14,-11) \hbox{ or }(14,-10,0,12,-11). 
$$
\end{proposition}

\begin{proof}
We have $\Lambda_{11}=\sfc$ and $\Lambda_{16}=0$ so $\sfc=0$: also
$\Lambda_{66}=\sfd+2\sfe-\sfb$ so $\sfb=\sfd+2\sfe$. Moreover
$\sfa=\Lambda_{55}\ge 0$ and $\sfe=-\Lambda_{45}\le 0$. Putting
$\sfx=\sfa+\sfb$ and $\sfy=\sfa+\sfe$ we find $\Lambda_{23}=\sfx-\sfy$
and $\Lambda_{63}=\sfy$, and both are positive by Lemma~1.1(ii), so
$\sfx>\sfy>0\ge\sfe$. But the double-point formula gives
\begin{eqnarray}\label{dbpG4}
0&=&-2\sfy^2+2\sfy\sfx-2\sfe\sfx+\sfx^2+6\sfe-11\sfx\nonumber\\
&=&\sfx(\sfx-11)-\sfe(2\sfx-6)+2\sfy(\sfx-\sfy)
\end{eqnarray}
which implies $\sfx<11$. Since equation~(\ref{dbpG4}) has no integer
solutions with $\sfx=3$, $\sfe$ is also bounded and a simple search
gives the solutions claimed, along with the solution $(3,4,0,8,-2)$
which is excluded because $\Lambda_{33}$ is odd. 
\end{proof}

For $\cG_5$ we have $D_2\equiv D_3$ and $D_5\equiv D_6$, and we choose
the basis $\beta_1=D_3^2$, $\beta_2=D_3D_6$, $\beta_3=D_6^2$,
$\beta_4=D_6D_7$, $\beta_5=D_7^2$. 

\begin{proposition}\label{G5}
   The class of any totally nondegenerate abelian surface in a toric Fano
  $4$-fold of type $\cG_5$ satisfies $\sfc=0$, $\sfd>0$,
  $\sfb=\sfd+\sfx$, $\sfa\ge 2\sfx>0$ and
  $3\sfa=2\sfd\sfx+\sfd^2-4\sfx-9\sfd$, where $\sfx=\sfa+\sfe$. 
\end{proposition}
\begin{proof}
  $\Lambda_{15}=0$ so $\Lambda_{11}$ and $\Lambda_{55}$ both vanish:
  they are equal to $\sfc$ and $\sfd+\sfx-\sfb$ respectively. 
  $\sfd=\Lambda_{23}$ and $\sfx=\Lambda_{35}$ and both are the only
  non-zero values in some row of $\Lambda$. The double-point formula
  gives the equation for $\sfa$ in terms of $\sfd$ and $\sfx$, and the
  remaining inequality comes from the fact that
  $\Lambda_{33}=\sfa-2\sfx$. 
\end{proof}

Again it might be possible to analyse these cases further by
considering maps to toric surfaces. 

\section{Types $\cI_9$, $\cL_{11}$, $\cL_{13}$, $\cQ_{16}$, $\cU_8$,
  $\cV_4$, $\cW$, $\cZ_1$  and $\cZ_2$}\label{typeother}

In these miscellaneous cases the calculations do not usually have much
in common: even when they do, as in the cases of $\cZ_1$ and $\cZ_2$,
the outcomes can be very different. 

We start with the case of $\cI_9$, where the restrictions are
strong. In this case $\rk\Pic X=4$. 
The primitive relations are 
\begin{equation*}
\begin{array}[center]{lll}
x_7+x_8=x_3,&x_3+x_6=x_7,&x_6+x_8=0,\\
x_1+x_2=x_7,&x_3+x_4+x_5=x_8,&x_4+x_5+x_7=0.   
\end{array}  
\end{equation*}
We have $D_1\equiv D_2$ and $D_4\equiv D_5$. We choose the basis
$\beta_1=D_4D_6$, $\beta_2=D_4D_7$, $\beta_3=D_4D_8$,
$\beta_4=D_6D_7$, $\beta_5=D_7^2$, $\beta_6=D_8^2$. 

\begin{proposition}\label{I9}
  If $\sfa\beta_1+\cdots+\sff\beta_6$ is the class of a totally
  nondegenerate abelian surface in a toric Fano $4$-fold of
  type $\cI_9$ then $(\sfa,\sfb,\sfc, \sfd, \sfe,\sff)$ is one of
  $$
(18, 14, -12, -16, -7, -6), (16, 8, -6, -12, -4, -3), (26, 6,
  -4, -16, -3, -2),
$$
$$
(14, 10, -6, -12, -5, -3) \hbox{ or } (26, 6, -2, -16,
  -3, -1). 
$$ 
\end{proposition}

\begin{proof}
$\Lambda_{26}=0$ so $\Lambda_{66}=2\sff-\sfc$ vanishes by
Lemma~1.1(iii) and the rows $\Lambda_{2*}$ and $\Lambda_{6*}$ are
proportional by Lemma~1.1(v). Since
$\Lambda_{25}=-\sfb-2\sfe$ and $\Lambda_{65}=0$ this gives
$\sfb=-2\sfe$. Moreover $\Lambda_{55}=\sfe-\sfa-2\sfd$ vanishes. 

Now it is convenient to use the variables $\sff=-\Lambda_{61}$,
$\sfx=2\sfe-\sfd=\Lambda_{51}$ and $\sfy=\sff-\sfe=\Lambda_{21}$: as
these are the only non-zero values in those rows, we have $\sfx>0$,
$\sfy>0$ and $\sff<0$. The double-point formula gives
\begin{eqnarray}\label{dbpI9}
0&=&-2\sff\sfy+2\sfy\sfx+\sfy^2-2\sff\sfx-5\sfx+4\sff-7\sfy\nonumber\\
&=&\sfy(\sfy-7)+\sfx(2\sfy-5)+(-\sff)(2\sfy-4)+2(-\sff)\sfx
\end{eqnarray}
so $\sfy\le 6$. From this, by solving equation~(\ref{dbpI9}) for
$\sff$, it is trivial to obtain the result claimed. 
\end{proof}

For $\cL_{11}$ the restrictions are rather weak. In this case $\rk\Pic
X=3$. The primitive relations are
\begin{equation*}
x_1+x_8=0,\ x_2+x_3=0,\ x_4+x_5=x_3,\ x_6+x_7=x_2.   
\end{equation*}
We have $D_1\equiv D_8$, $D_4\equiv D_5$ and $D_6\equiv D_7$. We
choose the basis $\beta_1=D_3D_5$, $\beta_2=D_3D_7$, $\beta_3=D_3D_8$,
$\beta_4=D_5D_7$, $\beta_5=D_5D_8$, $\beta_6=D_7D_8$.

\begin{proposition}\label{L11}
  If $\sfa\beta_1+\cdots+\sff\beta_6$ is the class of a totally
  nondegenerate abelian surface in a toric Fano $4$-fold of
  type $\cL_{11}$ then $\sfa=\sfb>0$, $\sfc=0$, $\sfe=\sff>0$ and
  $\sfd=\sfb\sff-2\sfb-2\sff\ge 0$. 
\end{proposition}

\begin{proof}
Since $\Lambda_{12}=0$, we have $0=\Lambda_{11}=\sff-\sfe$ and
$0=\Lambda_{22}=\sfe-\sff-2\sfc$, so $\sfe=\sff$ and $\sfc=0$. The
rows $\Lambda_{1*}$ and $\Lambda_{2*}$ are proportional by
Lemma~1.1(v), and since $\Lambda_{13}=\Lambda_{23}=\sff$ we have
either $\sff=0$ or $\sff>0$ and $\Lambda_{15}=\Lambda_{25}$. But
$\Lambda_{15}=\sfd$ and $\Lambda_{25}=\sfa+\sfb-\sfd$, so if $\sff\neq
0$ we have $\sfa=\sfb$. 

The double-point formula gives
$\sfa\sff+\sfb\sff-3\sfa-\sfb-2\sfd-4\sff=0$. This immediately
excludes $\sff=0$ and hence gives the formula for $\sfd$. Since
$\sfb=\Lambda_{35}$ it is non-negative, and $\sfb\neq 0$ because
otherwise $\sfd<0$. 
\end{proof}

For $\cL_{13}$ the restrictions are much stronger. In this case $\rk\Pic
X=3$. The primitive relations are
\begin{equation*}
x_1+x_8=0,\ x_2+x_3=x_1,\ x_4+x_5=x_1,\ x_6+x_7=x_8.   
\end{equation*}
We have $D_2\equiv D_3$, $D_4\equiv D_5$ and $D_6\equiv D_7$ and we
choose the basis $\beta_1=D_3D_5$, $\beta_2=D_3D_7$, $\beta_3=D_5D_7$,
$\beta_4=D_5D_8$, $\beta_5=D_7D_8$, $\beta_6=D_8^2$. 

\begin{proposition}\label{L13}
  If $\sfa\beta_1+\cdots+\sff\beta_6$ is the class of a totally
  nondegenerate abelian surface in a toric Fano $4$-fold of
  type $\cL_{13}$ then $(\sfa,\sfb,\sfc,\sfd,\sfe,\sff) =
  (2,2,0,4,4,0)$, $(6,3,3,0,6,2)$ or $(2,0,2,-4,8,4)$. 
\end{proposition}
\begin{proof}
Since $\Lambda_{15}=0$ we have $0=\Lambda_{11}=\sfa-\sfb-\sfc$ and
also $\Lambda_{55}=0$, which gives $\sfe=\sfd+3\sff$. The double-point
formula gives
\begin{eqnarray}\label{dbpL13}
0&=&2\sfb\sfd+3\sfb\sff+\sfc\sfd+3\sfc\sff+\sfd^2
+3\sfd\sff+3\sff^2-4\sfb-4\sfc-6\sfd-12\sff\nonumber\\
&=&2\sfb\sfx+(\sff-4)\sfb+\sfc\sfx+(2\sff-4)\sfc+\sfx^2+(\sff-6)\sfx+(\sff-6)\sff
\end{eqnarray}
where $\sfx=\sfd+\sff=\Lambda_{24}$. Since $\sfc=\Lambda_{12}$,
$\sfb=\Lambda_{13}$ and $\sff=\Lambda_{34}$ and they are thus
nonnegative, we deduce that $\sff\le 6$. Moreover $\sfb$ and $\sfc$
are not both zero, since otherwise $\Lambda_{1j}=0$ for all $j$
because $\Lambda_{14}=\sfb+\sfc$. 

By Lemma~1.1(v) we have
$\Lambda_{12}\Lambda_{53}=\Lambda_{13}\Lambda_{52}$, and as
$\Lambda_{25}=\sfc+\sff$ and $\Lambda_{53}=\sfb+\sfx$ we have
$\sfb\sff=\sfc\sfx$. This, equation~(\ref{dbpL13}) and the bound
$0\le\sff\le 6$ is sufficient to give the result by a simple
calculation. 
\end{proof}

For $\cQ_{16}$ we get strong restrictions. In this case $\rk\Pic
X=5$. The primitive relations are
\begin{equation*}
  \begin{array}[center]{llll}
x_8+x_9=0,&x_7+x_9=x_1,&x_1+x_2=x_9,&x_1+x_8=x_7,\\
x_2+x_7=0,&x_3+x_5=x_9,&x_4+x_6=x_8.&\\     
  \end{array}
\end{equation*}
We have $D_3\equiv D_5$, $D_4\equiv
D_6$. We choose the basis $\beta_1=D_5D_6$, $\beta_2=D_5D_7$,
$\beta_3=D_5D_8$, $\beta_4=D_5D_9$, $\beta_5=D_6D_7$,
$\beta_6=D_6D_8$, $\beta_7=D_6D_9$, $\beta_8=D_7^2$. 

\begin{proposition}\label{Q16}
  If $\sfa\beta_1+\cdots+\sfh\beta_8$ is the class of a totally
  nondegenerate abelian surface in a toric Fano $4$-fold of
  type $\cQ_{16}$ then $(\sfa,\ldots\sfh)$ is
  one of $(0,
  6, 7, -4, 6, 7, -4, 0)$, $(0, 4, 6, -2, 4, 6, -2, 0)$ or $(0, 3, 7,
  -1, 3, 7, -1, 0)$. 
\end{proposition}
\begin{proof}
  The computation shows that $\Lambda_{57}=\Lambda_{67}=0$ so by
  Lemma~1.1(iii) we have
  $\Lambda_{55}=\Lambda_{66}=\Lambda_{77}=0$. These give
  $-\sfa+\sfb-\sfe+2\sfh = -\sfb+\sfe = -\sfa-2\sfd+2\sfg =0$. Using
  these equations to eliminate $\sfa$, $\sfe$ and $\sfg$ we are left
  with $\Lambda_{17}=-\Lambda_{34}=\sfh$, so $\sfh=0$. Also
  $\Lambda_{12}=0$, and therefore $\Lambda_{1*}$ and $\Lambda_{2*}$
  are proportional; but $\Lambda_{16}=0$ and $\Lambda_{26}=\sff-\sfc$,
  so $\sff=\sfc$. Now it is convenient to put
  $\sfx=\sfb+\sfd=\Lambda_{13}>0$ (it is not zero, by Lemma~1.1(iv))
  and $\sfy=\sfc-\sfb=\Lambda_{53}>0$ similarly. We also have
  $0<\Lambda_{73}=-\sfd$. The double-point formula gives
\begin{eqnarray}\label{dbpQ16}
0&=&\sfx^2-2\sfx\sfz+2\sfy\sfx-2\sfy\sfd-6\sfx+4\sfd-4\sfy\nonumber\\
&=&(\sfx-6)\sfx+(2\sfx-4)(-\sfd)+(2sfx-4)\sfy+2\sfy(-\sfd)
\end{eqnarray}
so $1\le\sfx\le 5$. By solving equation~(\ref{dbpQ16}) for $\sfy$ or
$\sfd$ it is easy to see that the only solutions are as claimed. 
\end{proof}

For $\cU_8$ we again get strong restrictions. In this case $\rk\Pic
X=6$. The primitive relations are
\begin{equation*}
  \begin{array}[center]{llll}
x_1+x_3=x_2,&x_2+x_4=x_3,&x_1+x_4=0,&x_3+x_5=x_4,\\
x_4+x_6=x_5,&x_2+x_5=0,&x_1+x_5=x_6,&x_2+x_6=x_1,\\
x_3+x_6=0,&x_7+x_8=x_1,&x_9+x_{10}=x_4.& 
\end{array}
\end{equation*}
We have $D_7\equiv D_8$, $D_9\equiv D_{10}$. We choose the basis
$\beta_1=D_3D_8$, $\beta_2=D_3D_{10}$, $\beta_3=D_4D_8$,
$\beta_4=D_4D_{10}$, $\beta_5=D_5D_8$, $\beta_6=D_5D_{10}$,
$\beta_7=D_6^2$, $\beta_8=D_6D_8$, $\beta_9=D_6D_{10}$,
$\beta_{10}=D_8D_{10}$. 
\begin{proposition}\label{U8}
  If $\sfa\beta_1+\cdots+\sfj\beta_{10}$ is the class of a totally
  nondegenerate abelian surface in a toric Fano $4$-fold of
  type $\cU_8$ then $(\sfa,\sfx,\sfy)= (4,2,1)$ or $(2,2,2)$ up to
  permutation, and $(\sfa,\ldots,\sfj)=(\sfa, \sfa, \sfa+\sfx,
  \sfa+\sfx, \sfy+\sfx, \sfy+\sfx, 0, \sfy, \sfy, 0)$
\end{proposition}
\begin{proof}
  The computation shows that $\Lambda_{13}=0$ so by
  Lemma~1.1(iii) we have
  $\Lambda_{11}=\Lambda_{33}=0$. These give $\sfi=\sfh-\sfj$ and
  $\sfb=\sfa-\sfj$. The proportionality between $\Lambda_{1*}$ and
  $\Lambda{_3*}$ together with $\Lambda_{16}=0$ and
  $\Lambda_{36}=\sfj$ gives $\sfj=0$. Also $\Lambda_{14}=0$ so
  $\Lambda_{44}=0$ and by the proportionality between $\Lambda_{1*}$
  and $\Lambda_{4*}$ we have $\Lambda_{43}=0$. These give $\sfd=\sfc$
  and $\sfe=\sff$. Finally, $\Lambda_{55}=2\sfg=-\Lambda_{56}$ so
  $\sfg=0$. 

Now we write the double-point formula, in terms of
$\sfa=\Lambda_{72}$, $\sfx=\sff-\sfh=\Lambda_{76}$,
$\sfy=\sfa-\sfc+\sff=\Lambda_{74}$ and $\sfz=\sfc-\sff$. By
Lemma~1.1(iv) we get $\sfa, \sfx, \sfy>0$. The double-point formula
gives
\begin{equation}\label{dbpU8}
\sfz^2+(2\sfa+2\sfx+2\sfy-6)\sfz+(2\sfa\sfx+2\sfx\sfy+2\sfy\sfa-4(\sfa+\sfx+\sfy))=0
\end{equation}
which is symmetric in $\sfa$, $\sfx$ and $\sfy$. So to solve it we may
assume that $\sfa\ge\sfx\ge\sfy\ge -\sfz$ (since
$\sfy+\sfz=\Lambda_{81}>0$) and $\sfy>0$. It is straightforward to
check that the only solutions of equation~(\ref{dbpU8}) satisfying
these conditions are those given. 
\end{proof}

For $\cV_4$ we get some interesting restrictions. In this case
$\rk\Pic X=6$. The primitive relations are
\begin{equation*}
\begin{array}[center]{c}
x_4+x_{10}=0,\ x_1+x_5=0,\ x_2+x_6=0,\ x_3+x_7=0,\ x_8+x_9=0,\\
\begin{array}[center]{ll}
x_1+x_2+x_{10}=x_7+x_8,&x_1+x_3+x_{10}=x_6+x_8,\\
x_2+x_3+x_{10}=x_5+x_8,&x_1+x_2+x_3=x_4+x_8,\\
x_1+x_9+x_{10}=x_6+x_7,&x_2+x_9+x_{10}=x_5+x_7,\\
x_3+x_9+x_{10}=x_5+x_6,&x_1+x_2+x_9=x_4+x_7,\\
x_1+x_3+x_9=x_4+x_6,&x_2+x_3+x_9=x_4+x_5,\\
x_4+x_5+x_6=x_3+x_9,&x_4+x_5+x_7=x_2+x_9,\\
x_4+x_6+x_7=x_1+x_9,&x_5+x_6+x_7=x_9+x_{10},\\
x_4+x_5+x_8=x_2+x_3,&x_4+x_6+x_8=x_1+x_3,\\
x_4+x_7+x_8=x_1+x_2,&x_5+x_6+x_8=x_3+x_{10},\\
x_5+x_7+x_8=x_2+x_{10},&x_6+x_7+x_8=x_1+x_{10}. 
\end{array}\\
\end{array}  
\end{equation*}
No two $D_i$ are numerically
equivalent, so $\Lambda=L$. We choose the basis $\beta_1=D_5D_6$,
$\beta_{2}=D_5D_7$, $\beta_{3}=D_5D_8$, $\beta_{4}=D_5D_9$,
$\beta_{5}=D_5D_{10}$, $\beta_{6}=D_6D_7$, $\beta_{7}=D_6D_8$,
$\beta_{8}=D_6D_9$, $\beta_{9}=D_6D_{10}$, $\beta_{10}=D_7D_8$,
$\beta_{11}=D_7D_9$, $\beta_{12}=D_7D_{10}$, $\beta_{13}=D_8^2$,
$\beta_{14}=D_9^2$, $\beta_{15}=D_9D_{10}$, $\beta_{16}=D_{10}^2$. 
Thus $L$ is a $10\times 10$ matrix of linear forms in $16$ variables
$\sfa,\ldots,\sfp$. 
\begin{proposition}\label{V4}
  If $\sfa\beta_1+\cdots+\sfp\beta_{16}$ is the class of a totally
  nondegenerate abelian surface in a toric Fano $4$-fold of
  type $\cV_4$ then there is a symmetric matrix
\begin{equation*}
M=\begin{pmatrix}
0&\sfq&\sfu&\sfx&\sfz\\
\sfq&0&\sfr&\sfv&\sfy\\
\sfu&\sfr&0&\sfs&\sfw\\
\sfx&\sfv&\sfs&0&\sft\\
\sfz&\sfy&\sfw&\sft&0
\end{pmatrix}
\end{equation*}
with non-negative off-diagonal entries that determines
$\sfa,\ldots,\sfp$ by
\begin{eqnarray*}
(\sfa,\ldots,\sfp)&=&(\sfs+\sft+\sfw,\sft+\sfv+\sfy,
\sfr-\sft,\sfs+\sft+\sfv,
\sfr+\sfw+\sfy,\sft+\sfx+\sfz,\\
&&{}-\sft+\sfu,
\sfs+\sft+\sfx,\sfu+\sfw+\sfz,\sfq-\sft,
\sft+\sfv+\sfx,\sfq+\sfy+\sfz,\sft,\\
&&\sfs+\sft+\sfv+\sfx,
\sfq+\sfr+\sfu,\sfw+\sfy+\sfz). 
\end{eqnarray*}
The rank of $M$ is at most~$4$, no row of $M$ vanishes, and  
\begin{equation}\label{dbpV4}
\sum_{\{i,j\}\cap\{k,l\}=\emptyset} M_{ij}M_{kl}=4\sum_{1\le i,j\le 5} M_{ij}. 
\end{equation}
\end{proposition}
\begin{proof}
  There are off-diagonal entries of $L$ that vanish in each row, so
  every diagonal entry is zero by Lemma~1.1(iii). This gives a system
  of ten linear equations in the sixteen variables $\sfa,\ldots,\sfp$,
  which is of rank six, so ten of the sixteen variables are
  independent. The matrix $M$ is obtained by selecting the submatrix
  of $L$ given by the first four rows and columns and the eighth
  (which is independent of the first four for general values of
  $(\sfa,\ldots,\sfp)$). Then we set $\sfq=M_{1,2}$, etc., and use
  these ten variables subsequently. 

Equation~(\ref{dbpV4}) is the double-point formula and the condition
that $0<\rk M <5$ comes from Lemma~1.1(i). If a row of $M$
vanishes then (it is simple to check by computation) the corresponding
row of $L$ vanishes also, contrary to Lemma~1.1(iii). 
\end{proof}

For $\cW$ we get some restrictions. In this case $\rk\Pic X=5$. 
The primitive relations are
\begin{equation*}
  \begin{array}[center]{c}
  \begin{array}[center]{lll}
x_1+x_4=x_7,&x_2+x_5=x_8,&x_3+x_6=x_9,\\
x_1+x_2+x_3=0,&x_4+x_5+x_6=0,&x_7+x_8+x_9=0,\\
x_1+x_2+x_9=x_6,&x_4+x_5+x_9=x_3,&x_1+x_3+x_8=x_5,\\
x_4+x_6+x_8=x_2,&x_2+x_3+x_7=x_4,&x_5+x_6+x_7=x_1,
  \end{array}\\
  \begin{array}[center]{ll}
x_1+x_8+x_9=x_5+x_6,&x_4+x_8+x_9=x_2+x_3,\\
x_2+x_7+x_9=x_4+x_6,&x_5+x_7+x_9=x_1+x_3,\\
x_3+x_7+x_8=x_4+x_5,&x_6+x_7+x_8=x_1+x_2.    
  \end{array}
  \end{array}
\end{equation*}
No two $D_i$ are equivalent, so
$\Lambda=L$.  We choose the basis $\beta_1=D_3^2$, $\beta_2=D_3D_9$,
$\beta_3=D_6^2$, $\beta_4=D_6D_7$, $\beta_5=D_6D_8$, $\beta_6=D_6D_9$,
$\beta_7=D_7^2$, $\beta_8=D_7D_8$, $\beta_9=D_7D_9$,
$\beta_{10}=D_8^2$, $\beta_{11}=D_8D_9$, $\beta_{12}=D_9^2$. Thus $L$
is a $9\times 9$ matrix of linear forms in $12$ variables
$\sfa,\ldots,\sfl$. 

\begin{proposition}\label{W}
  If $\sfa\beta_1+\cdots+\sfl\beta_{12}$ is the class of a totally
  nondegenerate abelian surface in a toric Fano $4$-fold of
  type $\cW$ then there are positive integers $\sfw$, $\sfx$, $\sfy$ and
  an integer $\sft>-\min(\sfw,\sfx,\sfy)$ such that
\begin{equation}\label{dbpW}
3\sfh=\sft^2+(2\sfw+2\sfy+2\sfx-9)\sft+2\sfy(\sfw+\sfx)+2\sfw\sfx-7(\sfw+\sfx)-4\sfy. 
\end{equation}
Moreover $\sfh\ge 0$ and $\sfh\pm(\sfw-\sfx)$, $\sfh\pm(\sfx-\sfy)$,
$\sfh\pm(\sfy-\sfw)$ and $\sfh+\sfw+\sfx-2\sfy$ are also all non-negative. 

These numbers determine $(\sfa,\ldots,\sfl)$ by
\begin{eqnarray*}
(\sfa,\ldots,\sfl)&=&(\sfh+2\sft+\sfx+\sfw, 2\sft+\sfw+\sfx+\sfy,
\sfh+\sfx+\sfw, \sft,\\
&& \sft, -\sft+\sfw+\sfx+\sfy,
-\sfh-\sft-\sfw, \sfh,\\
&& \sfh+\sft+\sfw-\sfy, -\sfh-\sft-\sfx, \sfh+\sft+\sfx-\sfy, -\sfh+2\sfy). 
\end{eqnarray*}
\end{proposition}

\begin{proof}
The entries $L_{14}$, $L_{25}$ and $L_{36}$ vanish, so Lemma~1.1(iii) gives six
linear relations $L_{11}=\cdots=L_{66}=0$. Computing $L_{ii}$ we find
\begin{eqnarray*}
-\sff+\sfg+\sfj+\sfl&=&0\\
\sfa-\sfb+2\sfc+\sfe-\sff+\sfg+\sfj+\sfl&=&0\\
3\sfa-2\sfb+\sfg+\sfj+\sfl&=&0\\
2\sfa-\sfb+\sfc-\sfd-\sff+\sfg+\sfj+\sfl&=&0\\
2\sfa-\sfb+\sfc-\sfe-\sff+\sfg+\sfj+\sfl&=&0\\
3\sfc-2\sff+\sfg+\sfj+\sfl&=&0. 
\end{eqnarray*}
This system has rank~$4$ and we use it to
eliminate $\sfc$, $\sfe$, $\sff$ and $\sfg$. 

Note that $\sfd\neq 0$, because if $\sfd=0$ then $L_{58}=\sfa+\sfj$ and
$L_{88}=-\sfa-2\sfj$ so $\sfa=\sfj=0$, and then 
$$
L_{8*}=(-\sfh,0,-\sfk,-\sfh,0,-\sfk,\sfh,0,\sfk)
$$
so $\sfh=\sfk=0$ and $L_{8*}$ vanishes. 

Because $L_{14}=L_{25}=0$ we have, from Lemma~1.1(v), that
$L_{19}L_{48}=L_{18}L_{49}$, $L_{17}L_{48}=L_{18}L_{47}$ and
$L_{29}L_{58}=L_{28}L_{59}$. After dividing by $\sfd$ these give three
independent linear equations
\begin{eqnarray*}
-\sfh-\sfj-\sfb +\sfi+\sfl&=&0\\
2\sfa - 2\sfb +\sfl -\sfh&=&0\\
\sfa-\sfb+\sfj+\sfk+\sfl&=&0. 
\end{eqnarray*}
We use them to eliminate $\sfa$, $\sfb$ and $\sfl$. Now we introduce
$\sfw=L_{17}$, $\sfx=L_{27}$ and $\sfy=L_{37}$, which are all positive
(otherwise a row of $L$ vanishes) and
$\sft=L_{47}-\sfw=L_{57}-\sfx=L_{67}-\sfy$. We have $L_{78}=\sfh$. The
other non-zero values occurring are the ones listed in the theorem as
being non-negative, and equation~(\ref{dbpW}) is the double-point
formula. 
\end{proof}

For $\cZ_1$ we get strong restrictions. In this case $\rk\Pic
X=4$. The primitive relations are
\begin{equation*}
\begin{array}[center]{lll}
x_1+x_2+x_5=0,&x_1+x_2+x_6=x_7,&x_2+x_4+x_5=x_8,\\
x_2+x_4+x_6=x_7+x_8,&x_3+x_7+x_8=0,&x_3+x_4+x_6=x_1+x_5,\\
x_3+x_4+x_7=x_1,&x_3+x_6+x_8=x_5.&
\end{array}
\end{equation*}
No two $D_i$ are equivalent, so
$\Lambda=L$. We choose the basis $\beta_1=D_5^2$, $\beta_2=D_5D_6$,
$\beta_3=D_6^2$, $\beta_4=D_6D_7$, $\beta_5=D_6D_8$, $\beta_6=D_7^2$,
$\beta_7=D_7D_8$, $\beta_8=D_8^2$. Thus $L$ is an $8\times 8$ matrix
of linear forms in $8$ variables $\sfa,\ldots,\sfh$. 

\begin{proposition}\label{Z1}
  If $\sfa\beta_1+\cdots+\sfh\beta_8$ is the class of a totally
  nondegenerate abelian surface in a toric Fano $4$-fold of
  type $\cZ_1$ then $(\sfw,\sfx,\sfy,\sfz)=(0,2,2,6)$, $(0,3,1,8)$,
  $(1,2,1,2)$, $(0,1,3,8)$ or $(1,1,2,2)$. These integers determine
  $(\sfa,\ldots,\sfh)$ by
\begin{eqnarray*}
(\sfa,\ldots,\sfh)&=&
(\sfw+\sfx+\sfy+\sfz, 3\sfw+2\sfx+2\sfy+\sfz, 3\sfw+\sfx+\sfy,\\
&& 3\sfw+2\sfx, 2\sfy+\sfz, \sfw+\sfx, 2\sfy+2\sfz, -\sfy-\sfz). 
\end{eqnarray*}
\end{proposition}

\begin{proof}
The entries $L_{18}$ and $L_{57}$ vanish, so Lemma~1.1(iii) gives four
linear relations $L_{11}=L_{55}=L_{77}=L_{88}=0$, namely
\begin{eqnarray*}
\sfa-\sfb+\sfc-\sfd+2\sff&=&0\\
3\sfa-2\sfb+\sfc+\sfe+2\sfh&=&0\\
\sfc-2\sfd-\sfe+3\sff+\sfg+\sfh&=&0\\
2\sfa-\sfb+\sfc-\sfd+\sff+\sfg+3\sfh&=&0. 
\end{eqnarray*}
This system has rank~$3$ and we use it to
eliminate $\sfa$, $\sfb$ and $\sfc$. 

Since $L_{18}=0$ we have $L_{14}L_{85}=L_{15}L_{84}$ by 
Lemma~1.1(v) and this gives $\sfg=-2\sfh$. Now we may use
the variables $\sfw=L_{15}$, $\sfx=L_{16}$, $\sfy=L_{45}$ and
$\sfz=L_{46}$. The double-point formula gives
\begin{equation}\label{dbpZ1}
(2(\sfx+\sfy)+6\sfw-6)\sfz+(3\sfw-17)\sfw
+(6\sfw-11)(\sfx+\sfy)+(\sfx+\sfy)^2+4\sfx\sfy=0
\end{equation}
from which it follows immediately that $\sfw\le 5$. Then it is simple
to list the integer solutions for $(\sfw,\sfx,\sfy,\sfz)$. Note that
$L$ and equation~(\ref{dbpZ1}), but not $\sfa,\ldots,\sfh$, are
symmetric in $\sfx$ and $\sfy$ so we may assume that $\sfx\ge \sfy$ at
this stage. The solutions $(0,6,0,5)$, $(0,11,0,0)$, $(1,1,0,9)$ and
$(1,1,1,4)$ all give odd values for some $L_{ii}$, and $(0,4,0,14)$
is excluded because $L_{5*}$ vanishes. $(1,7,0,0)$ gives $L_{44}<0$
and $(1,2,0,5)$ gives $L_{45}=0$ but $L_{44}>0$, so both are excluded. 

The remaining solutions are $(0,2,2,6)$, $(0,3,1,8)$ and $(1,2,1,2)$
which, after possibly interchanging $\sfx$ and $\sfy$, give the result claimed. 
\end{proof}

For $\cZ_2$ we get only weak restrictions. In this case $\rk\Pic X=4$. 
The primitive relations are
\begin{equation*}
  \begin{array}[center]{lll}
x_1+x_2+x_5=0,&x_1+x_2+x_6=x_7,&x_2+x_4+x_5=x_8,\\
x_2+x_4+x_6=x_7+x_8,&x_3+x_7+x_8=x_2,&x_3+x_4+x_6=0,\\
x_3+x_4+x_7=x_1+x_2&x_3+x_6+x_8=x_2+x_5&\\
  \end{array}  
\end{equation*}
No two $D_i$ are equivalent, so
$\Lambda=L$. We choose the basis $\beta_1=D_5^2$, $\beta_2=D_5D_6$,
$\beta_3=D_6^2$, $\beta_4=D_6D_7$, $\beta_5=D_6D_8$, $\beta_6=D_7^2$,
$\beta_7=D_7D_8$, $\beta_8=D_8^2$. Thus $L$ is an $8\times 8$ matrix
of linear forms in $8$ variables $\sfa,\ldots,\sfh$. 

\begin{proposition}\label{Z2}
  If $\sfa\beta_1+\cdots+\sfh\beta_8$ is the class of a totally
  nondegenerate abelian surface in a toric Fano $4$-fold of
  type $\cZ_2$ then there are positive integers $\sfx$, $\sfy$, $\sfz$ and
  a non-negative integer $\sfw$ such that
\begin{equation}\label{dbpZ2}
3\sfw=\sfx^2+2\sfx\sfy+2\sfy\sfz+2\sfz\sfx-9\sfx-4\sfy-10\sfz. 
\end{equation}
Moreover $\sfw+\sfz\ge \sfy$ and $\sfw+2\sfz\ge 2\sfy$. 

These numbers determine $(\sfa,\ldots,\sfh)$ by
$$
(\sfa,\ldots,\sfh)=(\sfw+2\sfz, \sfx+\sfy+\sfz, 0, -\sfy, \sfx+\sfy,
0, \sfw+\sfx+\sfz, -\sfw-\sfz)
$$
\end{proposition}
\begin{proof}
  $L_{18}=0$ so $\sff=L_{11}=0$. But $L_{12}=\sff$ so $L_{22}=0$. Also
  $L_{57}=0$ so $L_{55}=L_{77}=0$. This gives two independent
  equations, which we use to eliminate $\sfa$ and $\sfd$. Applying
  Lemma~1.1(v) to the first two rows gives $\sfc=0$, since
  $L_{15}=\sfc$. Now we use the variables $\sfx=L_{32}$, $\sfy=L_{37}$
  and $\sfz=L_{38}$, all of which are positive (otherwise a row
  vanishes), and $\sfw=L_{44}$. Equation~(\ref{dbpZ2}) is the
  double-point formula and the inequalities are imposed by
  $L_{46}=\sfw-\sfy+\sfz$ and $L_{46}=\sfw-2\sfy+2\sfz$. 
\end{proof}
\section{Summary}
We conclude with a table that shows what is known about the existence
of totally nondegenerate abelian surfaces in smooth toric Fano
$4$-folds. The Fano $4$-folds are listed in the order of the tables
in~\cite{Sat3}. In the second column, $\checkmark$ indicates that it is known
that such surfaces do exist and $\cross$ indicates that it is known that
they do not. The symbols $\flat$ and $\sharp$ indicate that it is not
known whether such surfaces exist: $\sharp$ is used when there are
known to be only finitely many classes that could possibly accommodate
such a surface. The third column in each block gives a reference to a
paper or to a theorem in this paper, where more details may be found.
\begin{center}
\begin{tabular}{|c|c|c||c|c|c||c|c|c|}
\hline
%\multicolumn{9}{|c|}
%{Table of results}\\
%\hline
$\PP^4$&$\checkmark$&\cite{HM}&$\cG_1$&$\cross$&\cite{Sat2}&$\cQ_{10},\cQ_{11}$&$\checkmark$&\cite{Sat2}\\
\hline
$\cB_1$-$\cB_3$&$\cross$&\cite{Sat2}&$\cG_2$&$\cross$&\cite{Kaj2}&$\cQ_{12}-\cQ_{15}$&$\cross$&\cite{Sat2}\\
\hline
$\cB_4$&$\checkmark$&\cite{Lan}&$\cG_3$&$\flat$&\ref{G3}&$\cQ_{16}$&$\sharp$&\ref{Q16}\\
\hline
$\cB_5$&$\cross$&\cite{San}&$\cG_4$&$\sharp$&\ref{G4}&$\cQ_{17}$&$\cross$&\cite{Sat2}\\
\hline
$\cC_1,\cC_2$&$\cross$&\cite{San}&$\cG_5$&$\flat$&\ref{G5}&$\cK_1-\cK_3$&$\cross$&\cite{Sat2}\\
\hline
$\cC_3$&$\cross$&\ref{C3}&$\cG_6$&$\cross$&\cite{Kaj2}&$\cK_4$&$\checkmark$&\cite{Sat2}\\
\hline
$\cC_4$&$\checkmark$&\cite{Sat2}&$\cH_1-\cH_7$&$\cross$&\cite{Sat2}&$\cR_1-\cR_3$&$\cross$&\cite{Sat2}\\
\hline
$\cE_1-\cE_3$&$\cross$&\cite{Sat2}&$\cH_8$&$\checkmark$&\cite{Sat2}&$\cP$&$\cross$&\cite{Sat2}\\
\hline
$\cD_1-\cD_6$&$\cross$&\cite{Sat2}&$\cH_9,\cH_{10}$&$\cross$&\cite{Sat2}&$\cU_1-\cU_4$&$\cross$&\cite{Sat2}\\
\hline
$\cD_7$&$\sharp$&\ref{D7}&$\cL_1-\cL_6$&$\cross$&\cite{Sat2}&$\cU_5$&$\checkmark$&\cite{Sat2}\\
\hline
$\cD_8,\cD_9$&$\cross$&\cite{Sat2}&$\cL_7-\cL_9$&$\checkmark$&\cite{Sat2}&$\cU_6,\cU_7$&$\cross$&\cite{Sat2}\\
\hline
$\cD_{10}$&$\cross$&\ref{D10}&$\cL_{10}$&$\cross$&\cite{Sat2}&$\cU_8$&$\sharp$&\ref{U8}\\
\hline
$\cD_{11}$&$\flat$&\ref{D11}&$\cL_{11}$&$\flat$&\ref{L11}&$\Tilde\cV_4$&$\cross$&\cite{Kaj2}\\
\hline
$\cD_{12}$&$\cross$&\cite{Sat2}&$\cL_{12}$&$\cross$&\cite{Sat2}&$\cV_4$&$\flat$&\ref{V4}\\
\hline
$\cD_{13}$&$\checkmark$&\cite{Sat2}&$\cL_{13}$&$\sharp$&\ref{L13}&$S_2\times S_2$&$\checkmark$&\cite{Sat2}\\
\hline
$\cD_{14}$&$\flat$&\ref{D14}&$\cI_1-\cI_8$&$\cross$&\cite{Sat2}&$S_2\times S_3$&$\checkmark$&\cite{Sat2}\\
\hline
$\cD_{15}$&$\checkmark$&\cite{Sat2}&$\cI_9$&$\sharp$&\ref{I9}&$S_3\times S_3$&$\checkmark$&\cite{Sat2}\\
\hline
$\cD_{16}$&$\cross$&\cite{Sat2}&$\cI_{10}-\cI_{15}$&$\cross$&\cite{Sat2}&$\cZ_1$&$\sharp$&\ref{Z1}\\
\hline
$\cD_{17}$&$\flat$&\ref{D17}&$\cM_1-\cM_5$&$\cross$&\cite{Sat2}&$\cZ_2$&$\flat$&\ref{Z2}\\
\hline
$\cD_{18}$&$\cross$&\ref{D18}&$\cJ_1.\cJ_2$&$\cross$&\cite{Sat2}&$\cW$&$\flat$&\ref{W}\\
\hline
$\cD_{19}$&$\cross$&\cite{Kaj2}&$\cQ_1-\cQ_9$&$\cross$&\cite{Sat2}&&&\\
\hline
\end{tabular}
\end{center}

\bigskip

\noindent
G.K.~Sankaran,\\
Department of Mathematical Sciences,\\
University of Bath,\\
Bath BA2 7AY,\\
England\\
{\tt gks@maths.bath.ac.uk}
\end{document}